\documentclass[11pt,a4paper]{article}
\usepackage{amsmath}
\usepackage{amssymb}
\usepackage{amsthm}
\usepackage{amsfonts}
\usepackage{latexsym}
\usepackage{graphicx}
\usepackage{color}
\usepackage{amscd}
\usepackage{epsfig}
\usepackage{cite}

\newtheorem{problem}{Problem}[section]

\newtheorem{lemma}[problem]{Lemma}
\newtheorem{theorem}[problem]{Theorem}
\newtheorem{prop}[problem]{Proposition}
\newtheorem{corollary}[problem]{Corollary}

\newtheorem{remark}{Remark}

\numberwithin{equation}{section}
\allowdisplaybreaks[4]
\begin{document}
\title{$(1^3+1)(2^3+1)\cdots(n^3+1)$ is not a cube\footnote{Supported by National Natural Science Foundation of China (Grant No. 11401285)}}
\author{Chuanze Niu\footnote{School of
Mathematical Sciences, Liaocheng University, Liaocheng 252059,
P. R. China, E-mail: niuchuanze@lcu.edu.cn}}
\date{}
\maketitle
\begin{abstract}
For a positive integer $n,$ define
$$C_n=\prod_{k=1}^n(k^3+1).$$
In this paper we prove that there are no cubes in the integer sequence $C_n,~n=1,2,\cdots.$
\end{abstract}



\section{Introduction}
Let $f(x)=ax^d+b$ be a polynomial of degree $d$ in $\mathbb{Z}[x].$
Define the consecutive product of the sequence $f(k),~k=1,2,\cdots,$ by $C_n(f)=\prod_{k=1}^nf(k).$ It is interesting to discuss whether
$C_n(f)$ is a $d$-th power of integers or not.

When $f(x)=x^2+1,$ Javier Cilleruelo \cite{Jc} proved that $C_n(f)$ is a square only for $n=3$ using the
analytic methods and answered the conjecture proposed by \cite{T}.
More generally, Yang, et al. \cite{YT} showed that there are finitely many squares in $C_n(f)$ when $d=2$ and $(a,b)=1,$ $1\leq a\leq 10$ and $1\leq b\leq 20.$
Under $ABC$ conjecture, Zhang \cite{Zh} proved $C_n(f)$ is not a $d$-th power for $n$ large enough, where $d=2^l\cdot3^m$ with $l\geq 1, l+m\geq 2, ab\neq 0$ and such that $f(k)\neq 0$ for any $k\geq 1$. When $f(x)=x^2-1,$ Hong \cite{Hong} proved that there are infinitely
many $n$ such that $C_n(f)$ is a square using Dirichlet's unit
theorem.

As can be seen, when $f$ is of degree 2, there exists many celebrate results.
In this paper,  we assume $f(x)=x^3+1$ and let $C_n=C_n(f),$ drawing on the idea of Cilleruelo \cite{Jc}, we show that $C_n$ is never a cube, the main result is the following theorem.
\begin{theorem}\label{mainth}
For any positive integer $n,$ the value $C_n=(1^3+1)(2^3+1)\cdots(n^3+1)$ is
not a cube.
\end{theorem}
As a corollary of the main theorem, we get the following result.
\begin{corollary}
The Diophantine equation $$y^3=\prod_{k=1}^n(k^3+a^3)$$
has no solutions with $a=1.$
\end{corollary}

\section{Preliminaries}

In this section, we do some preparations that used to prove the main theorem of this paper which include the upper bound
of the $p$-adic order of $C_n$ and the upper bound of the distribution function
$\pi(x)$ of primes. We always assume $n\geq4$ since one may check $C_1,~C_2$ and $C_3$ are not cubes easily.

Let $p$ be a rational prime number, the $p$-adic order of a positive integer $n$ is denoted by $\text{ord}_p(n),$
then
$n=\prod_{}p^{\text{ord}_p(n)},$ where the product runs through all prime factors $p$ of $n.$ Note that
$$\text{ord}_p(n)=\sum_{1\leq j\leq\frac{\log n}{\log p}}\chi(p^j|n),$$
where $\chi(p^j|n)$ equals 1 or 0 depending on whether $n$ is divisible by $p^j$ or not, then

\begin{align}\label{compute equation of the p-adic order of C_n}
ord_{p}(C_n)=\sum_{1\leq j\leq\frac{\log(n^3+1)}{\log p}}\sharp\{1\leq k\leq n|~p^j|(k^3+1)\}.
\end{align}

\begin{lemma}\label{order of prime=n+1}
If $p=n+1$ is a prime factor of $C_n$, then $\text{ord}_p(C_n)=1$.
\end{lemma}
\proof
Since $p=n+1$, we have $ord_p(n^3+1)\geq1.$ If $ord_p(n^3+1)>1$ then $p|(n^2-n+1)$ and therefore $p$ divides
$(n+1)^2-(n^2-n+1)=3n,$ which is a contradiction since $p>n>3.$
We claim that there are no integers $k$ with $1\leq k<n$ such that $p|((n-k)^3+1).$ Otherwise, $p$ divides $n^3+1-(n-k)^3-1=3nk(n-k),$ which is impossible.
Therefore,
$$ord_{p}(C_n)=ord_p(n^3+1)=1.$$
This completes the proof.\endproof
In the proof of Lemma \ref{order of prime=n+1}, we see in the sequence $k^3+1,~k=1,2,\cdots,$ the first term that divisible by a prime $p$ is $(p-1)^3+1.$

\begin{lemma}\label{order of prime>n+1}
If $p>n+1$ is a prime factor of $C_n,$ then $\text{ord}_p(C_n)\leq2.$
\end{lemma}
\proof
Note that $p>k+1$ and $ p^2>k^2-k+1$, we see $p^2$ can not divide $k^3+1$ for positive integers $k\leq n.$
Therefore $\text{ord}_p(k^3+1)=0$ or 1.

If $p|(k^3+1)$ and   $p|(l^3+1)$ for $1\leq l<k\leq n,$ then $p|(k^2-k+1)$ and $p|(l^2-l+1)$ which yields that $p|(k-l)(k+l-1).$
Therefore $p|(k+l-1)$ and it leads to $k+l-1=p$ since $1\leq k+l-1<2p.$ So $k$ and $l$ are determined by each other. The number of those $k$ such that $p|(k^3+1)$
is at most 2, thus by Equation (\ref{compute equation of the p-adic order of C_n}),
$\text{ord}_p(C_n)\leq 2.$ The proof is done.
\endproof
\begin{prop}\label{upper bound of p}
Let $p$ be a prime factor of a cube $C_n,$ then $p\leq n$.
\end{prop}
\proof
Note that for a prime factor $p$ of a cube $C_n$, $\text{ord}_p(C_n)\geq3.$ The proposition follows from Lemmas \ref{order of prime=n+1} and \ref{order of prime>n+1}.
\endproof
By Proposition \ref{upper bound of p}, one can rewrite a cube $C_n$ as the product of primes $p\leq n$, i.e.,
$$C_n=\prod_{p\leq n}p^{\text{ord}_p(C_n)}.$$

Next we always assume $f(x)=x^3+1$ and we shall discuss the number of the roots of $f(x)\equiv 0(\text{mod~} p^{j})$ for any positive integer $j.$
\begin{lemma}\label{pneq3}
$f(x)(\text{mod~} p)$ has 3 different roots if $p\equiv 1(\text{mod}~3)$ while it has exactly 1 root if $p\equiv 2(\text{mod}~3).$
\end{lemma}
\proof
Note that $x\equiv -1(\text{mod~} p)$ is  a trivial solution of $f(x)\equiv 0(\text{mod~} p)$ and it is the unique solution if $p=2.$ For $p\geq5$, consider the congruence
\begin{align}\label{equation of degree 2}
x^2-x+1\equiv0(\text{mod~} p).
\end{align}
It is equivalent to
\begin{align}
(2x-1)^2\equiv-3(\text{mod~} p).
\end{align}
Therefore Equation (\ref{equation of degree 2}) has no  solutions if the Lengerde symbol $\big(\frac{-3}{p}\big)=-1$, while it has 2
different solutions if $\big(\frac{-3}{p}\big)=1$.
By Gauss reciprocity law,
$$\Big(\frac{-3}{p}\Big)=\Big(\frac{-1}{p}\Big)\Big(\frac{3}{p}\Big)=(-1)^{\frac{p-1}{2}}(-1)^{\frac{p-1}{2}}\Big(\frac{p}{3}\Big)=\Big(\frac{p}{3}\Big).$$
The lemma follows since $\big(\frac{p}{3}\big)$ equals $-1$ if $p \equiv 2(\text{mod}~3)$ and equals 1 if  $p \equiv 1(\text{mod}~3)$.
\endproof

\begin{lemma}\label{modulo power of p}
If $p\neq3,$  then the number of the roots of $f(x) (\text{mod~} p)$ equals to the number of the roots of $f(x)(\text{mod~} p^{j})$ for any positive integer $j.$
\end{lemma}
\proof
Note that $f(x)=x^3+1\equiv0(\text{mod~}p)$ and $f'(x)=3x^2\equiv0(\text{mod~}p)$  have no common solutions if $p\neq3$, then the lemma follows by Hensel' Lemma (or see\cite{HW}, section 8.3).
\endproof
\begin{lemma}\label{p=3}
$x^3+1\equiv 0(\text{mod}~3^{j})$ has 3 different roots for any integer $j\geq2.$
\end{lemma}
\proof
Work by induction on $j.$

When $j=2,$ the roots of $x^3+1\equiv 0(\text{mod}~3^2)$ are $x\equiv2,~5,~8(\text{mod~}3^2).$

Assume the roots of $x^3+1\equiv 0(\text{mod}~3^{j})$ are $x\equiv3^{j-1}-1,~2\cdot3^{j-1}-1,~3^{j}-1(\text{mod~}3^{j}).$
Next we show the roots of $x^3+1\equiv 0(\text{mod}~3^{j+1})$ are $x\equiv3^{j}-1,~2\cdot3^{j}-1,~3^{j+1}-1(\text{mod~}3^{j+1}).$ It suffices to show that $\xi^3+1\not\equiv0(\text{mod~}3^{j+1})$ when $\xi=3^{j-1}-1,~2\cdot3^{j-1}-1.$ Take $\xi=3^{j-1}-1$ as an example and the other case could be checked similarly. Note that $3j-3\geq j+1$ and $2j-1\geq j+1,$
$$(3^{j-1}-1)^3+1=3^{3j-3}-3^{2j-1}+3^{j}\equiv3^{j}\not\equiv0(\text{mod~}3^{j+1}).$$
This completes the proof of Lemma \ref{p=3}.
\endproof
\begin{remark}
In fact, we can also show for any prime $p$ and positive integer  $j,$ $x^p+1\equiv 0(\text{mod}~p^{j+1})$ has exactly $p$ different roots
$x\equiv ip^{j}-1(\text{mod~}p^{j+1}),~1\leq i\leq p.$
\end{remark}

\begin{prop}\label{number of the solutions}
Let $i$ be any integer. For any positive integer $j,$ $f(x)\equiv0(\text{mod~}p^{j})$ has only 1 solution in the interval $[i,i+p^{j}-1]$ if $p\equiv2(\text{mod~}3),$ whereas it has at most 3 different solutions in $[i,i+p^{j}-1]$ for other $p$'s.
\end{prop}
\proof
The proposition follows from Lemmas \ref{pneq3}, \ref{modulo power of p} and \ref{p=3}.
\endproof

\begin{remark}\label{consective divisiblity}
In the sequence $k^3+1,~k=1,~2,~\cdots,$ when $p\equiv2(\text{mod~}3)$, the first three terms that divisible by $p$ are $(p-1)^3+1,$ $(2p-1)^3+1$ and $(3p-1)^3+1.$ Therefore
$C_n$ is not a cube if $p-1\leq n\leq 3p-2.$
\end{remark}
As an immediate consequence of Equation (\ref{compute equation of the p-adic order of C_n}) and Proposition \ref{number of the solutions}, we get the upper bound of $\text{ord}_p(C_n).$
\begin{prop}\label{order of C_n}
We have
$$
\text{ord}_p(C_n)\leq\left\{
                   \begin{array}{ll}
                     \sum\limits_{1\leq j\leq\frac{\log(n^3+1)}{\log
p}}\big\lceil\frac{n}{p^j}\big\rceil, & \hbox{}p\equiv2(\text{mod~}3); \\
                     3\sum\limits_{1\leq j\leq\frac{\log(n^3+1)}{\log
p}}\big\lceil\frac{n}{p^j}\big\rceil, & \hbox{}p\not\equiv2(\text{mod~}3),
                   \end{array}
                 \right.
$$
where $\lceil x\rceil$ denotes the least integer that is greater than or equals to $x.$
\end{prop}

The inequality $C_n>(n!)^3$ which is the technique that is used in \cite{Jc}, that is, combining Proposition \ref{order of C_n}, we can deduce a bound for the sum of $\frac{\log p}{p-1}$ for prime numbers
$p\leq n$ in arithmetical progression when $C_n$ is a cube.
We end this section by recalling a result on the distribution of prime numbers. Let $\pi(x)$ be the number of primes which do not exceed $x.$
\begin{lemma}\label{upper bound of pi(n)}
If $x>10000,$ then
$$\pi(x)<\frac{1.139x}{\log x}.$$
\end{lemma}
\proof
See Lemma 2.5 in \cite{YT}.
\endproof

\section{Proof of Theorem \ref{mainth}}

It is easy to see that
$$\text{ord}_p(n!)=\sum_{1\leq j\leq\frac{\log n}{\log p}}\Big[\frac{n}{p^j}\Big],$$
where $[x]$ is the greatest integer not exceeding $x.$
When $p\leq n,$ we have
\begin{equation}\label{order of n!}
\text{ord}_p(n!)\geq\sum_{1\leq j\leq\frac{\log n}{\log p}}\Big(\frac{n}{p^j}-1\Big)\leq\frac{n-p}{p-1}-\frac{\log n}{\log p}\leq\frac{n-1}{p-1}-\frac{2\log n}{\log p}.
\end{equation}

{\it Proof of Theorem \ref{mainth}.}  Note that $$\prod_{p\leq n}p^{\text{ord}_p(C_n)}>\prod_{p\leq n}p^{3\text{ord}_p(n!)},$$
take logarithm on both sides, we get
$$\sum_{p\leq n}\text{ord}_p(C_n)\log p>3\sum_{p\leq n}\text{ord}_p(n!)\log p.$$
Then
$$\sum_{p\leq n\atop p\equiv2(\text{mod~3})}(\text{ord}_p(n!)-\frac{1}{3}\text{ord}_p(C_n))\log p<\sum_{p\leq n\atop p\not\equiv2(\text{mod~3})}(\frac{1}{3}\text{ord}_p(C_n)-\text{ord}_p(n!))\log p.$$

By Proposition \ref{order of C_n} and Equation (\ref{order of n!}),
when $p\not\equiv2(\text{mod~3}),$ we have
\begin{align*}
\frac{1}{3}\text{ord}_p(C_n)-\text{ord}_p(n!)&\leq\sum_{1\leq j\leq\frac{\log(n^3+1)}{\log p}}\Big\lceil\frac{n}{p^j}\Big\rceil-\sum_{1\leq j\leq\frac{\log n}{\log p}}\Big[\frac{n}{p^j}\Big]\\
&\leq\sum_{1\leq j\leq\frac{\log n}{\log p}}\Big(\Big\lceil\frac{n}{p^j}\Big\rceil-\Big[\frac{n}{p^j}\Big]\Big)+\sum_{\frac{\log n}{\log p}<j\leq\frac{\log (n^3+1)}{\log p}}\Big\lceil\frac{n}{p^j}\Big\rceil\\
&\leq\sum_{1\leq j\leq\frac{\log n}{\log p}}1+\sum_{\frac{\log n}{\log p}<j\leq\frac{\log (n^3+1)}{\log p}}1\\
&\leq \frac{\log (n^3+1)}{\log p};
\end{align*}
when $p\equiv2(\text{mod~3}),$ we have
\begin{align*}
\text{ord}_p(n!)-\frac{1}{3}\text{ord}_p(C_n)&\geq\sum_{1\leq j\leq\frac{\log n}{\log p}}\Big[\frac{n}{p^j}\Big]-\frac{1}{3}\sum_{1\leq j\leq\frac{\log(n^3+1)}{\log p}}\Big\lceil\frac{n}{p^j}\Big\rceil\\
&=\frac{2}{3}\sum_{1\leq j\leq\frac{\log n}{\log p}}\Big[\frac{n}{p^j}\Big]+\frac{1}{3}\sum_{1\leq j\leq\frac{\log n}{\log p}}(\Big[\frac{n}{p^j}\Big]-\Big\lceil\frac{n}{p^j}\Big\rceil)-\frac{1}{3}\sum_{\frac{\log n}{\log p}<j\leq\frac{\log(n^3+1)}{\log p}}\Big\lceil\frac{n}{p^j}\Big\rceil\\
&\geq\frac{2}{3}\cdot\frac{n-1}{p-1}-\frac{4}{3}\cdot\frac{\log n}{\log p}-\frac{1}{3}\cdot\frac{\log (n^3+1)}{\log p}.
\end{align*}
Therefore, we have
\begin{align*}
\frac{2(n-1)}{3}\sum_{p\leq n\atop p\equiv2(\text{mod~3})}\frac{\log p}{p-1}&\leq\sum_{p\leq n\atop p\not\equiv2(\text{mod~3})}\log(n^3+1)+\frac{4}{3}\sum_{p\leq n\atop p\equiv2(\text{mod~3})}\log n+\frac{1}{3}\sum_{p\leq n\atop p\equiv2(\text{mod~3})}\log (n^3+1)\\
&\leq\frac{2}{3}\sum_{p\leq n\atop p\not\equiv2(\text{mod~3})}\log(n^3+1)+\frac{4}{9}\sum_{p\leq n\atop p\equiv2(\text{mod~3})}\log (n^3+1)+\frac{1}{3}\sum_{p\leq n}\log(n^3+1)\\
&=\frac{2}{9}\sum_{p\leq n\atop p\not\equiv2(\text{mod~3})}\log(n^3+1)+\frac{4}{9}\sum_{p\leq n}\log (n^3+1)+\frac{1}{3}\sum_{p\leq n}\log(n^3+1)\\
&\leq\pi(n) \log(n^3+1).
\end{align*}
Therefore
$$\sum_{p\leq n\atop p\equiv2(\text{mod~3})}\frac{\log p}{p-1}\leq \frac{3}{2(n-1)}\pi(n) \log(n^3+1).$$

Now using the upper bound of $\pi(x)$ in Lemma \ref{upper bound of pi(n)}, when $n>10000,$
\begin{align*}
\sum_{p\leq n\atop p\equiv2(\text{mod~}3)}\frac{\log
p}{p-1}&\leq \frac{3.417}{2}\cdot\frac{n}{n-1}\cdot\frac{\log(n^3+1)}{\log n}.
\end{align*}
After computations, the right-hand side is $<5.2$ if $n\geq10000.$
Adding over enough prime numbers for $p\equiv2(\text{mod~}3),$ we see that for
$n\geq34631,$
$$\sum_{p\leq n\atop p\equiv2(\text{mod~}3)}\frac{\log
p}{p-1}\geq5.2.$$
Hence $C_n$ is not a cube when $n\geq34631.$ For $4\leq n\leq34631,$ by Remark \ref{consective divisiblity}, we consider  the following primes that congruent to 2 modulo 3,
$$5,~11,~29,~83,~233,~683,~2039,~6113,~18329,$$
which corresponding intervals such that $C_n$ are not cubes if $n$ lies in them:
$$[4,13],~[10,31],~[28,85],~[83,247],~[232,697],~[682,2047],~[2038,6115],~[6112,18337],~[18328,54985].$$
Therefor $C_n$ is not a cube if $4\leq n\leq54985,$ this completes the proof.


\begin{thebibliography}{}
\bibitem{Jc} J.Cilleruelo, Squares in $(1^2+1)\cdots(n^2+1)$,
J. Number Theory 128, 2008, 2488-2491.
\bibitem{T} T. Amdeberhan, et al., Arithmetical propeties of a sequence arising from an arctangent sum, J. Number Theory 128(6), 2008, 1807-1846.
\bibitem{YT} Yang S, et al., Diophantine equations with products of consecutive values of a quadratic polynomial, J. Number Theory. 2011, 131: 1840-1851.
\bibitem{Zh} Z.Zhang, Powers in $\prod\limits_{k=1}^n (ak^{2^l\cdot3^m}+b)$, Functiones Et Approximatio Commentarii Mathematici, 2012, 46(1): 7-13.
\bibitem{Hong} Shaofang Hong and Xingjiang Liu, Squares in
$(2^2-1)\cdots(n^2+1)$ and $p$-adic valuation, Asian-Europe J. Math, 2010, 3(2): 329.
\bibitem{HW} G. Hardy, E.Wright. An introducton to the theory of numbers, sixth edition, Oxford Univ. Press, 2009.
\end{thebibliography}
\end{document}